\documentclass{amsart}

\usepackage{amsmath}
\usepackage{graphicx}
\usepackage{amsfonts}
\usepackage{amssymb}
\usepackage[utf8]{inputenc}
\usepackage{fancyhdr}
\usepackage{amsaddr}
\usepackage{fancyvrb}
\usepackage[normalem]{ulem}

\usepackage{datetime}
\usepackage{amsthm}
\usepackage{algorithm, caption}
\usepackage[noend]{algpseudocode}
\usepackage{bbm}
\usepackage{xcolor}
\usepackage[hidelinks]{hyperref}
\usepackage{colonequals}
\usepackage{comment}
\usepackage{listings}
\usepackage{tikz}
\usepackage{tikz-cd}

\usepackage[all]{xy} \SelectTips{eu}{}

\newtheorem{thm}{Theorem}[section]

  \newtheorem{lem}[thm]{Lemma}

\newtheorem*{theorem*}{Theorem A}

  \theoremstyle{definition}

  \numberwithin{equation}{thm}

  \makeatletter
  \let\c@algorithm\c@thm
  \makeatother

\newcommand{\K}{\mathbbm{k}}

\newcommand{\m}{\mathfrak{m}}

\newcommand{\tagform}[1]{\textnormal{\small(#1)}}
\newcommand{\tagKFlin}{\tagform{KF1}}
\newcommand{\tagKFmax}{\tagform{KF2}}
\newcommand{\tagKFcol}{\tagform{KF3}}

\title{The Pinched Veronese has no Koszul Filtration}
\author{Emily Berghofer}\author{\vspace{-30pt}}
\address{Mälardalen University}
\email{emily.berghofer@mdu.se}

\date{\today}
\keywords{Pinched veronese, Koszul algebra, Koszul filtration}
\subjclass[2020]{Primary 13D02; Secondary 13F55, 16S37}

\begin{document}
\maketitle
\begin{abstract}
In this paper we prove that the pinched Veronese subalgebra\\ $\K [x^3,x^2y,xy^2,x^2z,xz^2,y^3,y^2z,yz^2,z^3] \subset \K[x,y,z]$ has no Koszul filtration regardless of the field $\K$, despite it being a Koszul algebra. We do this by showing that no ideal $(f,g)$ minimally generated by two linear forms can belong to a Koszul filtration. By the definition of a Koszul filtration this, in turn, implies that no Koszul filtration can exist.

\end{abstract}
\section{Introduction}
\noindent
Let $\K$ be a field. A Koszul algebra is a standard graded $\K$-algebra such that $\K$ has a linear minimal free resolution. The notion of Koszulness was first introduced by Priddy \cite{Priddy1970}. Given the infinite and unpredictable nature of a free resolution, it is quite difficult to verify Koszulness directly from its definition. One must instead use different methods. One method is to show that, for some presentation, the defining ideal of the algebra has a quadratic Gröbner basis. If this is the case, we call the algebra \emph{G-quadratic} \cite[Remark 2.6(4)]{conca2012koszulalgebrasregularity}. Another important method is to find a family of ideals within the algebra that forms a \emph{Koszul filtration} \cite{CTV2001}. Both G-quadraticity and the existence of a Koszul filtration are sufficient conditions for Koszulness, but neither is necessary \cite{CTV2001}.

One particularly elusive Koszul algebra is the cubic pinched Veronese algebra.
The pinched Veronese is defined as \[PV=\K [x^3,x^2y,xy^2,x^2z,xz^2,y^3,y^2z,yz^2,z^3]\subset \K[x,y,z]\] where $\K$ is a field. 
In 1993, Bernd Sturmfels asked Irena Peeva whether the pinched Veronese is Koszul and more than a decade later, Caviglia answered this question affirmatively \cite{Caviglia2009}. A more theoretical proof was subsequently given by Caviglia and Conca \cite{CC2013}. Both proofs had to go outside the more standard methods described above.

It is still unknown whether the pinched Veronese is G-quadratic, although one might conjecture the answer to be no. Using a Gröbner fan as first introduced by Mora and Robbiano \cite{mora_grobner_1988}, Conca, De Negri and Rossi showed that for toric coordinates it has no quadratic Gröbner basis \cite{conca2012koszulalgebrasregularity}. It is also unknown if it has a quadratic Gröbner basis after a lifting by some linear form, i.e. if it is \emph{LG-quadratic}. 

It was previously shown by this author along with Nicklasson, Thompson and Westerbäck that the pinched Veronese does not have a Koszul filtration by monomial ideals \cite{berghofer2026con}. It was previously unknown, however, whether it has a Koszul filtration in general. In this paper, we will add another layer to the mystery of the pinched Veronese as we show that despite it being Koszul, it does not have a Koszul filtration regardless of $\K$.

\begin{theorem*}\label{bigone}The pinched Veronese does not have a Koszul filtration for any~field~$\K$.
\end{theorem*}
\section{Preliminaries}
\noindent
Let $R=\bigoplus_{i\geq 0}R_i$ be a standard graded $\K$-algebra. A \emph{Koszul filtration} of $R$ is a family of ideals $\mathcal{F}$ satisfying:
\begin{enumerate}
    \item[\tagKFlin] Every ideal in $\mathcal{F}$ is generated by elements of $R_1$.
    \item[\tagKFmax] The maximal ideal $\m_R$ belongs to $\mathcal{F}$.
    \item[\tagKFcol] For every nonzero ideal $I\in \mathcal{F}$, there exists an ideal $J\in \mathcal{F}$ properly contained in $I$ such that $I/J$ is cyclic and $J:I\in \mathcal{F}$.
\end{enumerate}
Recall that for ideals $I,J \in R$, the colon ideal is $J:I=\{f\in R \mid fI \subseteq J\}$. We will use the notation $J:x$ for $J:(J+(x))$.

Note that by (KF3), $\mathcal{F}$ contains a chain of ideals from the zero ideal $(0)$ to the maximal ideal $\m_R$: 
$$(0)\subset I_1 \subset I_2\subset\dots\subset\m_R$$
\noindent
where $I_i$ is minimally generated by $i$ generators and such that each colon ideal $I_i:I_{i+1}$ is linear and again in $\mathcal{F}$. This will be important for the proof of Theorem A. For every ideal $I$ in a Koszul filtration, $R/I$ has a linear minimal free resolution \cite{conca2012koszulalgebrasregularity}.

Let $A=\K [x_1,...,x_n]$ be a polynomial algebra, and let $S\subset A$ be any subalgebra. The following Lemma describes the colon ideal $(f):g$ where $(f)$ and $(f,g)$ are both ideals in a subalgebra $S\subset A$. This description resembles other descriptions of colon ideals, such as Theorem $14$ in Section $4$, Chapter $4$ of Cox, Little, and O'Shea \cite{cox2015ideals}.

\begin{lem}\label{coli}
 Let $f$ and $g$ be two elements in the subalgebra $S\subset A$, where $f=f_0h$ and $g=g_0h$ such that $\gcd_A(f_0,g_0)=1$. We can describe the colon ideal $(f):_Sg$ as follows:
    $$(f):_Sg=\{f_0p\mid f_0p, g_0p \in S \text{ and } p\in A\}$$
\end{lem}
\begin{proof}
We first show that the containment $\subseteq$ holds. Recall that by the definition of a colon ideal, $a\in (f):_Sg$ implies that 
$$ag=bf \Longleftrightarrow ag_0h=bf_0h\Longleftrightarrow ag_0=bf_0$$

\noindent
for some $b\in S$. Since we are in a unique factorisation domain, this implies $f_0|a$ and $g_0|b$. Hence we can express $a=f_0p_1$ and $b=g_0p_2$.
$$f_0p_1g_0=g_0p_2f_0\Longleftrightarrow p_1=p_2$$
\noindent
Thus $a=f_0p_1$ belongs to the right-hand side. For the inverse containment $\supseteq$, take any 
$ c \in \{f_0p\mid f_0p, g_0p \in S \text{ and } p\in A\}$. Then $cg=(f_0p)(g_0h)=(f_0h)(g_0p) \in (f)$, hence $c\in (f):_Sg$.

\end{proof}

\section{Obstructions for a Koszul filtration}
\noindent
Let $\K$ be a field and let $PV=\K[x^3,x^2y,xy^2,x^2z,xz^2,y^3,y^2z,yz^2,z^3]\subset \K[x,y,z]$. Throughout this paper, \emph{PV-linear} refers to elements in $PV_3$; similarly, \emph{PV-quadratic} refers to elements in $PV_6$. We also refer to $A\colonequals\K[x,y,z]$. By abuse of notation we will talk about $PV$ having a Koszul filtration despite it not being standard graded. One should view this as a Koszul filtration of a $R$ such that $R_n=PV_{3n}$.

For a Koszul filtration $\mathcal{F}$ of $PV$ to exist, there must exist at least one ideal $I\in \mathcal{F}$ minimally generated by $PV$-linear forms $f,g$ such that the ideal $(f):g$ is $PV$-linear and such that $PV/(f,g)$ has a linear minimal free resolution as a $PV$-module. 
In this section, we show that no such ideal exists. 

\begin{lem}\label{asix}
For all $n\geq 2$ it is true that $PV_{3n}=A_{3n}$.
\end{lem}
\begin{proof}
We start by proving that $PV_{6}=A_{6}$. The following multiplication table shows that all $\binom{6+2}{2}=28$ degree $6$ monomials are present.
    \[
\resizebox{\textwidth}{!}{$
\renewcommand{\arraystretch}{1.4}
\begin{array}{c|ccccccccc}
 & x^3 & x^2y & x^2z & xy^2 & xz^2 & y^3 & y^2z & yz^2 & z^3 \\ \hline
x^3  & \uline{x^6} & \uline{x^5y} & \uline{x^5z} & x^4y^2 & x^4z^2 & x^3y^3 & x^3y^2z & x^3yz^2 & x^3z^3 \\
x^2y &     & \uline{x^4y^2} & \uline{x^4yz} & \uline{x^3y^3} & \uline{x^3z^3} & \uline{x^2y^4 }& x^2y^3z & x^2y^2z^2 & x^2yz^3 \\
x^2z &     &        & \uline{x^4z^2} & \uline{x^3y^2z} & \uline{x^3yz^2} & \uline{x^2y^3z} & x^2y^2z^2 & \uline{x^2yz^3} & x^2z^4 \\
xy^2 &     &        &        & \uline{x^2y^4} & \uline{x^2y^2z^2} & \uline{xy^5} & \uline{xy^4z} & xy^3z^2 & xy^2z^3 \\
xz^2 &     &        &        &        & \uline{x^2z^4} & \uline{xy^3z^2} & \uline{xy^2z^3} & \uline{xyz^4} & \uline{xz^5} \\
y^3  &     &        &        &        &        & \uline{y^6} & \uline{y^5z} & y^4z^2 & \uline{y^3z^3} \\
y^2z &     &        &        &        &        &      & \uline{y^4z^2} & y^3z^3 & y^2z^4 \\
yz^2 &     &        &        &        &        &      &         & \uline{y^2z^4} & \uline{yz^5} \\
z^3  &     &        &        &        &        &      &         &         & \uline{z^6}
\end{array}
$}
\]
\noindent
 We now prove that $PV_{3n}=A_{3n}$ by induction on $n$. The case $n=2$ was established above. Assume that
\[
A_{3(n-1)}=PV_{3(n-1)}
\]
for some $n\geq 3$. Clearly, $PV_{3n}\subseteq A_{3n}$ since $PV\subseteq A$. For the reverse inclusion, take any monomial $m\in A_{3n}$. Since $m$ has degree $3n$, we can factor it as
\[
m=m_1m_2,
\]
where $m_1\in A_6$ and $m_2\in A_{3(n-2)}$. Since $A_6=PV_6$, the degree-$6$ factor $m_1$ can be written as a product of two monomials in $PV_3$. Thus
\[
m\in PV_3\cdot\bigl(PV_3A_{3(n-2)}\bigr).
\]
Note that 
\[ PV_3A_{3(n-2)}\subseteq A_{3(n-1)}\] 
so
\[
m\in PV_3A_{3(n-1)}.
\]
By the above assumption $PV_3A_{3(n-1)}=PV_3PV_{3(n-1)}=PV_{3n}$ and hence
\[
A_{3n}\subseteq PV_{3n}.
\]
This gives
\[
A_{3n}=PV_{3n}.
\]
\end{proof}
\noindent
Throughout this section, we let $f$ and $g$ be two $PV$-linear forms that are not $\K$-linear multiples of each other, and define $f=f_0h$ and $g=g_0h$ where $\gcd_A (f,g)=h$, where $\gcd_A$ is defined in the polynomial ring $A$. Since $f$ and $g$ have degree $3$ in $A$, there are three possible cases for the degree of $h$. 
We begin with the case where $\deg\gcd_A(f,g)=0$.
\begin{lem}\label{deg0}
    If $\deg\gcd_A(f,g)=0$, the colon ideal $(f):g$ is not generated by PV-linear forms. In fact, it has a minimal $PV$-quadratic generator.
\end{lem}
\begin{proof}

Assume that $\deg\gcd_A(f,g)=0$. Without loss of generality, assume that $h=1$, this implies that $f=f_0$ and $g=g_0$. By Lemma \ref{coli} any $a\in (f):g$ must be divisible by $f$. Since $f_0$ is itself $PV$-linear, the $PV$-linear elements of the colon ideal are precisely the scalar multiples of $f$. 

The colon ideal will always contain $fxyz$. By Lemma \ref{asix} $PV_6=A_6$ and $\deg_A(f_0)=\deg_A(g_0)=3$ , hence $\deg_A(fxyz)=\deg_A(gxyz)=6$ which means that they are both in  $PV_6\subset PV$. 
Since $PV \subseteq A$ and $A$ is a domain, an equality $fxyz=fu$ with $u\in PV$ would imply $u=xyz$, contradicting $xyz\notin PV$. Hence the $(f):g$ cannot be generated by only $PV$-linear forms.
\end{proof}
\noindent
Assume instead that $\deg\gcd_A(f,g)=1$. This means that $\deg f_0 =\deg g_0 =2$

\begin{lem}\label{deg1}
If $\deg\gcd_A(f,g)=1$, the colon ideal $(f):g$ has only PV-linear generators if and only if $f_0,g_0$ are in $\text{span}\{x^2,y^2,z^2\}$.   
\end{lem}
\begin{proof}
Assume that $f_0,g_0$ are in $\text{span}\{x^2,y^2,z^2\}$. Then by  Lemma \ref{coli}:
$$(f):g=(f_0x,f_0y,f_0z)$$
\noindent
Conversely, suppose that at least one of $f_0$ and $g_0$
contains a mixed quadratic term. We show that $(f):g$ cannot be generated entirely by PV-linear forms.
If $f_0$ or $g_0$ contains any term in $A_2\backslash \text{span}\{x^2,y^2,z^2\}$, then for all $\ell\in \{x,y,z\}$, $f_0\ell^4,g_0\ell^4\in PV$ but for at least one $\ell\in{x,y,z}$, either $f_0\ell$ or $g_0\ell$ does not belong to $PV$. By Lemma \ref{asix} $PV_6=A_6$, but multiplying a mixed term by the third variable produces a multiple of $xyz$ which is not in $PV$.

$f_0z\notin PV_3$,  so $f_0z\notin(f):g$. On the other hand, $f_0z^4,g_0z^4\in A_6=PV_6$ , and hence $f_0z^4\in(f):g$. 

By Lemma \ref{coli}, any generator of $(f):g$ must be divisible by $f_0$. Hence, $f_0z^4\in(f):g$ cannot be  reduced to a $PV$-linear generator. The same argument applies to the mixed terms $xz$ and $yz$.
\end{proof}
\noindent

\begin{lem}\label{Q}
   If $\deg\gcd_A(f,g)=2$ there exists a nonzero $Q=axy+bxz+cyz$ such that both $f_0Q$ and $g_0Q$ are in $PV$.
\end{lem}

\begin{proof}
    Since $\deg(f_0)=\deg(g_0)= 1$ we can write
    \[f_0=\mu_1x+\mu_2y+\mu_3z\]
    \[g_0=\mu'_1x+\mu'_2y+\mu'_3z\]
    Since $f_0$ and $g_0$ are linearly independent, the following two equations have a one-dimensional solution space in $\K^3$. 
    \[
    \begin{cases}
      X\mu_3+Y\mu_2+Z\mu_1=0\\
      X\mu'_3+Y\mu'_2+Z\mu'_1=0
    \end{cases}
    \]
Choose a nonzero solution $(a,b,c)$ and set
          $$Q=a xy+b xz+c yz.$$
  
    Note that every other solution is a scalar multiple of $(a,b,c)$. Thus, the corresponding quadratic forms $Q$ are scalar multiples of one another. We now compute
        \[Qf_0=axy(\mu_1x+\mu_2y+\mu_3z)+bxz(\mu_1x+\mu_2y+\mu_3z)+cyz(\mu_1x+\mu_2y+\mu_3z)\]
and
        \[Qg_0=axy(\mu'_1x+\mu'_2y+\mu'_3z)+bxz(\mu'_1x+\mu'_2y+\mu'_3z)+cyz(\mu'_1x+\mu'_2y+\mu'_3z)\]
By the definition of $Q$, the $xyz$-terms vanish. Hence both $Qf_0$ and $Qg_0$ are in $PV$.
  
\end{proof}

\begin{lem}\label{deg2}
 Let $f$ and $g$ be $PV$-linear forms such that $\deg\gcd_A(f,g)=2$, and $(f):g$ is generated by $PV$-linear forms. Then $PV/(f,g)$ has a non-linear minimal free resolution.
 $$
\begin{tikzcd}
	{ }& F_3 & PV(-2)^4 & PV(-1)^2 & PV & {PV/(f,g)} & 0
	\arrow[from=1-1, to=1-2]
	\arrow["{\partial_3}", from=1-2, to=1-3]
	\arrow["{\partial_2}", from=1-3, to=1-4]
	\arrow["{\partial_1}", from=1-4, to=1-5]
    \arrow [from= 1-5, to=1-6]
    \arrow [from= 1-6, to=1-7]
\end{tikzcd}
$$
Where $F_3$ is a graded free module having a summand $PV(-4)$. In particular $\partial_3$ must contain $PV$-quadratic entries.
\end{lem}

\begin{proof}
We want to show that $\partial_3$ must contain $PV$-quadratic entries. We prove this by first describing $\partial_1$ and $\partial_2$ explicitly. We know that that one minimal choice for $\partial_1$ is
\[
\partial_1 =
\begin{pmatrix}
f & g
\end{pmatrix}.
\]

 Next we look at $\partial_2$, whose columns are  minimal generators for the module of syzygies of the rows of $\partial_1$. That is, the solutions of $\lambda_1f+\lambda_2g=0$ for $\lambda_i\in PV$. We note that $\lambda_1,\lambda_2$ are solutions to this equation only if $\lambda_2\in (f):g$. 

 By Lemma \ref{coli}, all solutions are of the form $(\lambda_1,\lambda_2)=(pg_0,-pf_0)$ such that $pg_0,pf_0\in PV$. The PV-linear elements of the colon ideal are those such that $p\in A_2$ and 
 \[
p\in(x^2,y^2,z^2,xy,xz,yz).
\]
If $p\in\operatorname{span}_k\{x^2,y^2,z^2\}$, then $pf_0,pg_0\in PV$. The proof of Lemma \ref{Q} shows that $pf_0$ and $pg_0$ are in $PV$ if and only if the mixed terms of $p$ are a scalar multiple of a nonzero solution $Q=axy+bxz+cyz$ where $a,b,c\in \K$. Hence the elements $p\in A_2$ for which both $pf_0$ and $pg_0$
belong to $PV$ are precisely the elements of
\[
\operatorname{span}_k\{x^2,y^2,z^2,Q\}.
\]
 Now we consider $p\in A_5$. Every
degree-five monomial belongs to
\[
PV_3x^2+PV_3y^2+PV_3z^2
\]
except
\[
x^3yz,\qquad xy^3z,\qquad xyz^3.
\]
We show that these belong to $\operatorname{span}_k\{x^2,y^2,z^2,Q\}$.

$$x^2zQ-bxz^2(x^2)-cx^2y(z^2)= ax^3yz$$
$$x^2yQ-axy^2(x^2)-cx^2z(y^2)= bx^3yz$$
$$x^3Q-bx^2z(x^2)-ax^2y(x^2)= cx^3yz$$

Since at least one of $a,b,c$ is nonzero, it follows that $x^3yz\in \operatorname{span}_k\{x^2,y^2,z^2,Q\}$.
By symmetry, the same holds for $xy^3z$ and $xyz^3$. Hence, $A_5=PV_3x^2+PV_3y^2+PV_3z^2+PV_3Q$.

If $n\geq 2$, every monomial in $A_{3n-1}$ is divisible by at least one of $x^2,y^2,z^2$. Hence, every element $f_0p\in PV_{3n}$ can be written as a sum of multiples of $f_0x^2,f_0y^2$, and $f_0z^2$, with the coefficients in $A_{3(n-1)}=PV_{3(n-1)}$. 

It follows that $(f):g=(f_0x^2,f_0y^2,f_0z^2,f_0Q)$,
and hence 

\[
\partial_2 =
\begin{pmatrix}
g_0x^2 & g_0y^2 & g_0z^2 & g_0Q\\
-f_0x^2 & -f_0y^2 & -f_0z^2 & -f_0Q
\end{pmatrix},
\]

We now look at $\partial_3$. 
The columns of $\partial_3$ minimally generate
of $\ker(\partial_2)$.  Thus, a column
\[
\begin{pmatrix}
\lambda_1\\
\lambda_2\\
\lambda_3\\
-\lambda_4
\end{pmatrix}
\]
for $\lambda_i\in PV$, belongs to $\ker(\partial_2)$ precisely when
\[
\begin{cases}
\lambda_1g_0x^2+\lambda_2g_0y^2+
\lambda_3g_0z^2-\lambda_4g_0Q=0,\\
\lambda_1f_0x^2+\lambda_2f_0y^2+
\lambda_3f_0z^2-\lambda_4f_0Q=0.
\end{cases}
\]
or equivalently
\begin{equation}\label{syz}
\lambda_1x^2+\lambda_{2}y^2+\lambda_{3}z^2
=\lambda_{4}Q.
\end{equation}
We show that every generating set for this syzygy module must contain minimal $PV$-quadratic generators.

Consider the following three solutions:
\[
C_x=
\begin{pmatrix}
Qx^4\\
0\\
0\\
-x^6
\end{pmatrix},
\qquad
C_y=
\begin{pmatrix}
0\\
Qy^4\\
0\\
-y^6
\end{pmatrix},
\qquad
C_z=
\begin{pmatrix}
0\\
0\\
Qz^4\\
-z^6
\end{pmatrix}.
\]
By Lemma~\ref{Q}, $Q\neq 0$. Moreover, all nonzero entries
of these columns are PV-quadratic by Lemma \ref{asix}.

We now show that at least one of these three syzygies cannot
be generated by syzygies whose entries are PV-linear.

Suppose that
\[
C_i=r_1D_1+\cdots+r_nD_n
\]
for some $i\in\{x,y,z\}$, where $r_j\in PV$ and each $D_j$
is a solution of \eqref{syz} whose entries are PV-linear. Write
\[
D_j=
\begin{pmatrix}
\mu_{1,j}\\
\mu_{2,j}\\
\mu_{3,j}\\
-\mu_{4,j}
\end{pmatrix},
\]
where each $\mu_{k,j}$ is $PV$-linear, and
\begin{equation}
\mu_{1,j}x^2+\mu_{2,j}y^2+\mu_{3,j}z^2
=\mu_{4,j}Q.
\end{equation}

Every $\mu_{k,j}$ can be written as a linear
combination of the nine PV-linear monomial generators. Hence
\[
\begin{aligned}
\mu_{4,j}={}&
\gamma_1x^3+\gamma_2x^2y+\gamma_3x^2z
+\gamma_4y^3+\gamma_5xy^2\\
&+\gamma_6y^2z+\gamma_7z^3
+\gamma_8xz^2+\gamma_9yz^2,
\end{aligned}
\]
where $\gamma_k\in\K$.

Since $xyz\notin PV$, none of the following three monomials $x^3yz$, $xy^3z$ and
$xyz^3$ are divisible by $x^2,y^2$ or $z^2$ such that the remainder is in $PV$. As we are in a UFD this means that the left-hand side of \eqref{syz} is unable to form terms with those three monomials. In particular, this means that the coefficients of these terms in $\mu_{4,j}Q$ must vanish.

Expanding $\mu_{4,j}Q$ gives us
$$
\begin{aligned} &\mu_{4,j}Q=\gamma_1x^3axy+\gamma_1x^3bxz+\gamma_1x^3yz+\gamma_2x^2yaxy+\gamma_2x^2ybxz+\gamma_2x^2ycyz\\
  &+\gamma_3x^2zaxy+\gamma_3x^2zbxz+\gamma_3x^2zcyz+\gamma_4y^3axy+\gamma_4y^3bxz+\gamma_4y^3cyz\\
  &+\gamma_5xy^2axy+\gamma_5xy^2bxz+\gamma_5xy^2cyz+\gamma_6y^2zaxy+\gamma_6y^2zbxz+\gamma_6y^2zcyz\\
  &+\gamma_7z^3axy+\gamma_7z^3bxz+\gamma_7z^3cyz+\gamma_8xz^2axy+\gamma_8xz^2bxz+\gamma_8xz^2cyz\\
  &+\gamma_9yz^2axy+\gamma_9yz^2bxz+\gamma_9yz^2cyz.
\end{aligned}
$$
The vanishing of terms containing $x^3yz$, $xy^3z$ and
$xyz^3$ is equivalent to
\begin{equation}\label{eq}
\begin{cases}
c\gamma_1+b\gamma_2+a\gamma_3=0,\\
c\gamma_5+b\gamma_4+a\gamma_6=0,\\
c\gamma_8+b\gamma_9+a\gamma_7=0.
\end{cases}
\end{equation}
We now evaluate $\mu_{4,j}$ at the point $(c,b,a)$. We have
\[
\mu_{4,j}(c,b,a)
=\gamma_1c^3+\gamma_2c^2b+\gamma_3c^2a
+\gamma_4b^3+\gamma_5cb^2+\gamma_6b^2a
+\gamma_7a^3+\gamma_8ca^2+\gamma_9ba^2
=\]\[c^2(\gamma_1c+\gamma_2b+\gamma_3a)
+b^2(\gamma_5c+\gamma_4b+\gamma_6a)
+a^2(\gamma_8c+\gamma_9b+\gamma_7a). 
\]
It follows from \eqref{eq} that
\[
\mu_{4,j}(c,b,a)=0.
\]
Thus, the bottom entry of every syzygy whose entries are
$PV$-linear vanishes at $(c,b,a)$.

Consequently, if $C_i$ were a $PV$-linear combination of
such syzygies, its bottom entry would also vanish at
$(c,b,a)$. However, the bottom entries of $C_x,C_y,C_z$
evaluate, respectively, to
\[
-c^6,\qquad -b^6,\qquad -a^6.
\]
Since $Q=axy+bxz+cyz\neq 0$, at least one of $a,b,c$ is
nonzero. Therefore, at least one of these three evaluations
is nonzero, giving a contradiction.

Since at least one of $C_x,C_y,C_z$ cannot be generated
by syzygies whose entries are $PV$-linear, the syzygy module $\ker \partial_2$ is not generated in $PV$-degree one. Hence the minimal free resolution of $PV/(f,g) $ is not linear.
\end{proof}

\section*{Proof of Theorem A}
\begin{proof}
    Assume that the pinched Veronese has a Koszul filtration $\mathcal{F}$.
    By (KF3) $\mathcal{F}$ must contain at least one ideal that is minimally generated by two $PV$-linear forms $f$ and $g$. By minimality, the two $PV$-linear forms are not scalar multiples of each other. 

    There are three potential scenarios for the relationship of $f$ and $g$ as cubic forms in $A$, either $\deg\gcd_A(f,g)=0,\:\deg\gcd_A(f,g)=1$ or $\deg\gcd_A(f,g)=2$.
    
     By Lemma \ref{deg0}, for all ideals $(f,g)$ such that $\deg\gcd_A(f,g)=0$ the colon ideal $(f):g$ cannot be generated by $PV$-linear forms. This violates (KF3) and hence no such ideal can belong to $\mathcal{F}$. By Lemma \ref{deg1}, the same conclusion holds when $\deg\gcd_A(f,g)=1$ and at least one of $g_0, f_0$ contains a nonzero term divisible by $xy$, $xz$, or $yz$.
     
      By Lemma $\ref{deg1}$, if $\deg\gcd_A(f,g)=1$ and both $f_0$ and $g_0$ belong to $\operatorname{span}_k\{x^2,y^2,z^2\}$, then $(f):g=(f_0x,f_0y,f_0z)$. Since $\deg f_0=2$, any subset ideal generated by two $PV$-linear forms $m$ and $n$ will have $\deg\gcd_A(m,n)=\deg(f_0)=2$.      
      
      By Lemma \ref{deg2} for all ideals $(f,g)$ such that $\deg\gcd_A(f,g)=2$, $PV/(f,g)$ has no linear free resolution. Therefore, $(f,g)$ is not in $\mathcal{F}$. 
      
      Thus $\mathcal{F}$ contains no ideal minimally generated by two elements, which is a contradiction. Hence no Koszul filtration of the pinched Veronese can exist.
\end{proof}

\noindent
The inspiration for the proof of Theorem A comes from exhaustive computations over finite fields. Using the algorithms from \cite{berghofer2026con}, it first became apparent that the algebra did not have any "obvious" Koszul filtration. Using Macaulay 2 \cite{M2}, all two-generated ideals, first over $\mathbb{F}_2$ and later over $\mathbb{F}_3$, were checked for the first few steps of their free resolutions. This was followed by an investigation of their possible colon ideals. The results followed the exact trend described for a general field $\K$ in Section 3 and hence inspired the proof. 

An example of a characteristic-dependent Koszul filtration was given in \cite{berghofer2026con}, but it seems to be due to the effect characteristic $2$ has on the defining ideal. An interesting open question is when there is a relationship between a Koszul filtration in characteristic $0$ and a Koszul filtration in other characteristics.

\section*{Acknowledgements}
\noindent
I would like to thank Lisa Nicklasson for introducing me to the Pinched Veronese and my supervisor Peder Thompson for his support and valuable discussions. 
The only use of artificial intelligence or large language models in this paper was during the editing process.

\bibliographystyle{amsplain}\bibliography{ref.bib}
\end{document}